\documentclass{amsart}
\usepackage{amssymb}

\theoremstyle{definition}

\theoremstyle{remark}

\numberwithin{equation}{section}

\input{tcilatex}

\begin{document}
\title{Complete Reduction in K-Theory}
\author{Mehmet K\i rdar}
\address{Department of Mathematics, Faculty of Arts and Science, Nam\i k
Kemal University, Tekirda\u{g}, Turkey}
\email{mkirdar@nku.edu.tr}
\date{June 8, 2011}

\begin{abstract}
We define a reduction, called complete reduction, for the K and KO relations
of the Hopf bundle over lens spaces introducing some numbers of interest to
various theories of mathematics. Along the way, we make an interesting
conjecture in number theory related to the cyclotomic fields.
\end{abstract}

\maketitle






%    General info

\section{Introduction}

Let $BZ_{n}$ denote the classifying space of the cyclic group $Z_{n}$. The $%
K $-ring of $BZ_{n}$ is given by $K(BZ_{n})=Z\left[ \mu \right] \diagup
((1+\mu )^{n}-1)$ where $\mu =\eta -1$ is the reduction of the tautological
complex line bundle (Hopf bundle) over $BZ_{n}.$ When $n$ is odd, the $KO$%
-ring of $BZ_{n}$ is described in the following way: Let $w=r(\mu )$ be the
realification of $\mu $. Then, $KO(BZ_{n})=Z\left[ \omega \right] \diagup
(wf_{n}(w))$ where $f_{n}(\cdot )$ is the polynomial 
\begin{equation*}
f_{n}(w)=n+\sum_{j=1}^{\frac{n-3}{2}}\frac{%
n(n^{2}-1^{2})(n^{2}-3^{2})...(n^{2}-(2j-1)^{2})}{2^{2j}.(2j+1)!}w^{j}+w^{%
\frac{n-1}{2}}.
\end{equation*}%
The topological $K$-theory over the real numbers is a little more involved
when $n$ is even. See [2] for details.

From now on, for simplicity, let $n=p$ be an odd prime number, although the
very simple idea of the paper can be extended for all natural numbers. In
this case, the topological K-theory of the lens spaces is very well-studied.
See [1].

In this note, we will define a reduction called "complete reduction" for the
relations of $\mu $ and $\omega $ coming from the generators of the
principal ideals of the above rings. Complete reduction is the smallest way
of writing these relations geometrically, respecting the Atiyah-Hirzebruch
spectral sequence and detects the group cohomology of $Z_{p}$ by means of
the filtrations of that spectral squence.

In order to obtain first few terms of the complete reduction, we make a
division trick and this gives some invariants -numbers- which we named $%
K_{n} $ for the complex case and $M_{n}$ for the real case, interesting not
only for $K\Lambda $-rings of lens spaces in topological $K$-theory but for $%
R\Lambda $-rings of cyclic groups in representation theory and for
cyclotomic rings of integers in number theory, due to the equivalence of
theories $K(BZ_{p}),$ $R(Z_{p})$ and $Z\left[ \exp \frac{2\pi i}{p}\right] $.

\section{K-Reduction}

By iteration, the relation $(1+\mu )^{p}-1=0$ can be put in the form

\begin{equation*}
p\mu =-\mu ^{p}+\frac{p-1}{2}\mu ^{p+1}+a_{2}\mu ^{p+1}+......+a_{n}\mu
^{p+n}+......
\end{equation*}%
\ \ \ \ \ \ \ \ \ \ \ \ \ \ \ 

\textrm{Definition 1.1. }The relation above is called completely reduced if $%
\left\vert a_{n}\right\vert \leq \frac{p-1}{2}$ for all $n\geq 2.$

Obviously the complete reduction is unique.

\textrm{Example 1.2.} For $p=3,$ the complete reduction is periodic with
period 2 and repeating coefficients $-1,1.$ For $p=5,$ the complete
reduction is periodic with period 6 and repeating coefficients $%
-1,2,-2,1,0,0.$ For $p=7,$ the first 28 coefficients of the complete
reduction is as below:%
\begin{eqnarray*}
7\mu &=&-\mu ^{7}+3\mu ^{8}+3\mu ^{9}+2\mu ^{10}+2\mu ^{11}+3\mu ^{12}+\mu
^{13}-2\mu ^{14}+0.\mu ^{15}+\mu ^{16} \\
&&+\mu ^{17}-2\mu ^{18}-2\mu ^{19}+0.\mu ^{20}+3\mu ^{21}-\mu ^{22}+2\mu
^{23}+\mu ^{24}-3\mu ^{25} \\
&&-\mu ^{26}+3\mu ^{27}+0.\mu ^{28}+2\mu ^{29}+0.\mu ^{30}+2\mu ^{31}-\mu
^{32}-2\mu ^{33}+\mu ^{34} \\
&&-653\mu ^{35}-3662\mu ^{36}-5800\mu ^{37}-4373\mu ^{38}-1651\mu
^{39}-253\mu ^{40}
\end{eqnarray*}%
Quite surprisingly, we couldn't observe periodicity of the coefficients of
the complete reduction for $p=7.$ One should do further reduction to decide
if it is periodic or not. Although we had the opposite belief after this
experimentation, as an interesting open problem in number theory, we make
the following probably false conjecture.

\textbf{Conjecture 1.3.}\textrm{\ }The complete reduction is periodic for
all odd prime numbers.

Next, we want to express first few coefficients of the complete reduction in
terms of $p.$ We introduce a very simple idea, probably done before, many
times in history. We will do the division trick used in the example above.

\textrm{Definition 1.4. }Define integers $K_{p,n}$ by 
\begin{equation*}
\sum_{n=0}^{\infty }K_{p,n}\mu ^{n}=\frac{-p\mu }{\left( 1+\mu \right)
^{p}-1-\mu ^{p}}
\end{equation*}

Let us denote $K_{p,n}$ simply by $K_{n}$ when $p$ is understood. Then $p\mu
=\sum_{n=0}^{\infty }K_{n}\mu ^{p+n}$ is a reduction of the relation of $\mu
.$ But, of course, it is not the complete reduction except for the primes $3$
and $5$. On the other hand, the first $p+1$ coefficients of the complete
reduction of the relation of $\mu $ are $K_{n}$ $(\func{mod}$ $p),$ $0\leq
n\leq p.$

The numbers $K_{n}$ satisfy a recursive formula. By using this recursive
formula or by direct division, which is the same process, we can compute $%
K_{n}$ for all $n\leq p-2$ as a polynomial of $p.$ We computed upto $K_{6}$
as below:

\begin{eqnarray*}
K_{1} &=&\frac{p-1}{2},\text{ \ \ for }p\geq 3 \\
K_{2} &=&-\frac{p^{2}-1}{12},\text{ \ \ for }p\geq 4 \\
K_{3} &=&\frac{p^{2}-1}{24},\text{ \ \ for }p\geq 5 \\
K_{4} &=&\frac{(p^{2}-1)(p^{2}-19)}{720},\text{ \ \ for }p\geq 6 \\
\text{ }K_{5} &=&-\frac{(p^{2}-1)(p^{2}-9)}{480},\text{ \ for }p\geq 7 \\
K_{6} &=&-\frac{(p-1)(2p^{5}+122p^{4}-1825p^{3}+8375p^{2}-17617p+15263)}{%
60480},\text{ \ \ for }p\geq 8
\end{eqnarray*}%
The author doesn't know whether these, probably very well-known, polynomials
are used somewhere in number theory. For large primes, we can use these
tabulated formulas, to find at least first $p-1$ terms of the complete
reduction for the prime number $p$.

\textrm{Example 1.5. }For $p=23,$ $K_{1}=11,$ $K_{2}=-44\equiv 2,$ $%
K_{3}=22\equiv -1,$ $K_{4}=374\equiv 6,$ $K_{5}=-572\equiv 3,$ $%
K_{6}=-10494\equiv -6,$ and hence the first seven terms of the complete
reduction are%
\begin{equation*}
23\mu =-\mu ^{23}+11\mu ^{24}+2\mu ^{25}-\mu ^{26}+6\mu ^{27}+3\mu
^{28}-6\mu ^{29}+......
\end{equation*}

Here, we recall the famous Bernoulli number $B_{n}$. It immediately follows
from the definitions that $B_{n}=\lim_{p\rightarrow \infty }\frac{-n!K_{n}}{%
p^{n}}.$

\section{KO-Reduction}

By iteration, the relation $wf_{p}(w)=0$, explicitly, 
\begin{equation*}
p\omega +\sum_{j=2}^{\frac{p-1}{2}}\frac{%
p(p^{2}-1^{2})(p^{2}-3^{2})...(p^{2}-(2j-3)^{2})}{2^{2j-2}.(2j-1)!}\omega
^{j}+\omega ^{\frac{p+1}{2}}=0
\end{equation*}%
can be written in the form%
\begin{equation*}
p\omega =-\omega ^{\frac{p+1}{2}}+b_{1}\omega ^{\frac{p+3}{2}}+b_{2}\omega ^{%
\frac{p+5}{2}}+......\text{ \ \ \ \ \ \ \ \ \ \ \ \ \ \ \ \ \ \ \ \ \ \ \ \
\ \ \ \ \ \ \ \ \ \ \ \ \ }
\end{equation*}

Similar to the complex case, we call the relation above completely reduced
if $\left\vert b_{n}\right\vert \leq \frac{p-1}{2}$ for all $n\geq 1.$
Complete reduction is clearly unique.

\textrm{Example 2.1. }For $p=3,$ the complete reduction is $3\omega =-\omega
^{2}$. It is periodic of period $1$, all coefficients, after $a_{0}=-1,$
being $0$. For $p=5,$ $5\omega =-\omega ^{3}+\omega ^{4}+5\omega ^{3}$ is
the complete reduction. It is periodic of period $2$ with repeating
coefficients $-1,+1.$ For $p=7,$ we did some reduction and found first $\ 16$
terms of the complete reduction as below:%
\begin{eqnarray*}
7\omega &=&-\omega ^{4}+2\omega ^{5}-3\omega ^{6}-3\omega ^{7}+2\omega
^{8}-\omega ^{9}-\omega ^{10}-3\omega ^{11} \\
&&-\omega ^{12}-\omega ^{13}+\omega ^{14}+\omega ^{15}+\omega ^{16}+\omega
^{17}-\omega ^{18}-3\omega ^{19} \\
&&-2481\omega ^{20}-1627\omega ^{21}-266\omega ^{22}
\end{eqnarray*}%
Again, similar to the complex case, we couldn't observe a periodicty for the
prime number $7.$ On the other hand, we conjecture that the complete
reduction is periodic in the real case too.

Next we define some numbers for the computation of the coefficients of the
complete reduction in terms of $p$.

\textrm{Definition 2.2. }Define integers $M_{p,n}$ by 
\begin{equation*}
\sum_{n=0}^{\infty }M_{p,n}\omega ^{n}=\frac{-p\omega }{wf_{p}(w)-\omega ^{%
\frac{p+1}{2}}}
\end{equation*}

Let us denote $M_{p,n}$ simply by $M_{n}$ when $p$ is understood. Then $%
p\omega =\sum_{n=0}^{\infty }M_{n}\omega ^{n+\frac{p+1}{2}}$ is a reduction
for $\omega .$ Of course, it is not complete reduction except for the primes 
$3$ and $5$. On the other hand, the first $\frac{p+1}{2}$ coefficients of
the complete reduction of the relation of $\omega $ are $M_{n}$ $(\func{mod}$
$p),0\leq n\leq \frac{p-1}{2}.$

Clearly $M_{0}=-1$ for all $p.$ We can calculate $M_{n}$ by writing a
recursive formula like we did in complex case, or by direct division. We
obtain formulas for $M_{n}$ in terms of $p$ which are valid for $p\geq 2n+3.$
The next three are%
\begin{eqnarray*}
M_{1} &=&\frac{p^{2}-1}{24},\text{ \ \ \ \ \ }p\geq 5 \\
M_{2} &=&-\frac{(p^{2}-1)(7p^{2}+17)}{5760},\text{ \ \ \ \ \ }p\geq 7 \\
M_{3} &=&\frac{(p^{2}-1)(57p^{4}-34p^{2}+169)}{322560},\text{ \ \ \ \ \ }%
p\geq 9
\end{eqnarray*}

\textrm{Example 2.3. }For $p=23,$ $M_{1}=22\equiv -1,$ $M_{2}=-341\equiv 4,$ 
$M_{3}=26081\equiv -1$ and the first four terms of the complete reduction are%
\begin{equation*}
23\omega =-\omega ^{12}-\omega ^{13}+4\omega ^{14}-\omega ^{15}+......
\end{equation*}

Note also that by means of the realification map, $M_{n}$ should be
expressed in terms of $K_{n}.$ So, if the conjecture is true for the complex
case, it is also true for the real case.

\end{document}